\documentclass[10pt,reqno]{amsart}
\usepackage{amsmath,amsfonts, amssymb, amsthm}
  \usepackage{paralist}
  \usepackage{graphics} 
  \usepackage{epsfig} 
\usepackage{graphicx}  \usepackage{epstopdf}
 \usepackage[colorlinks=true]{hyperref}
\hypersetup{urlcolor=blue, citecolor=red}

\newtheorem{theorem}{Theorem}[section]
\newtheorem{corollary}{Corollary}

\newtheorem{lemma}[theorem]{Lemma}

\theoremstyle{definition}

\newtheorem{remark}{Remark}

\def\pmod #1{\ ({\rm{mod}}\ #1)}
\def\Z{\Bbb Z}
\def\N{\Bbb N}

\def\R{\Bbb R}

\def\l{\left}
\def\r{\right}
\def\bg{\bigg}
\def\({\bg(}
\def\){\bg)}
\def\t{\text}
\def\f{\frac}

\def\Aut{\text{Aut}}

\def\gen{{\rm gen}}

\def\ls{\leq}
\def\gs{\geq}
\def\se {\subseteq}

\def\bi{\binom}
\def\al{\alpha}

\def\ve{\varepsilon}

\def\eq{\equiv}

\def\da{\delta}

\issueinfo{27}
{}
{}
{2019}
\PII{eISSN: 2688-1594 AIMS (2019)}
\dateposted{December 21, 2019} 
\pagespan{1}{19}

\copyrightinfo{2019}{American Institute of Mathematical Sciences}

\begin{document}
 \hbox{Electron. Res. Arch. 27 (2019), 69--87.}
\medskip
 
\title[Some universal quadratic sums over the integers] 
      {Some universal quadratic sums over the integers}

\author[Hai-Liang Wu and Zhi-Wei Sun]{Hai-Liang Wu and Zhi-Wei Sun$^*$}

\thanks{The initial version of this paper was posted to arXiv (with the code {\tt arXiv:1707.06223})
in July 2017. The second author is supported by the National Natural Science
Foundation of China (grant 11571162) and the NSFC-RFBR Cooperation and Exchange Program (grant 11811530072).}

\thanks{$^*$ Corresponding author: Zhi-Wei Sun}
\address{
Hai-Liang Wu,
Department of Mathematics, Nanjing University, Nanjing 210093, China}
 \email{whl.math@smail.nju.edu.cn}

\address{
Zhi-Wei Sun,
Department of Mathematics, Nanjing University, Nanjing 210093, China}
 \email{zwsun@nju.edu.cn}

%
%
%
%
%

\begin{abstract}
Let $a,b,c,d,e,f\in\mathbb N$ with $a\gs c\gs e>0$, $b\ls a$ and $b\equiv a\pmod2$, $d\ls c$ and $d\equiv c\pmod2$,
$f\ls e$ and $f\equiv e\pmod2$. If any nonnegative integer can be written as
$x(ax+b)/2+y(cy+d)/2+z(ez+f)/2$ with $x,y,z\in\mathbb Z$, then the ordered tuple $(a,b,c,d,e,f)$ is said to be universal over $\Z$.
Recently, Z.-W. Sun found all candidates for such universal tuples over $\Z$.
In this paper, we use the theory of ternary quadratic forms to show that 44 concrete tuples $(a,b,c,d,e,f)$ in Sun's list of candidates are indeed universal over $\mathbb Z$.
For example, we prove the universality of $(16,4,2,0,1,1)$ over $\Z$ which is related to the form $x^2+y^2+32z^2$.\vspace*{-20pt}
\end{abstract}
\maketitle

\section{Introduction}
Those $p_3(x)=x(x+1)/2$ with $x\in\Z$ are called triangular numbers.
In 1796 Gauss proved Fermat's assertion that each $n\in\N=\{0,1,2,\ldots\}$ can be expressed as the sum of three triangular numbers.

For polynomials $f_1(x),f_2(x),f_3(x)$ with $f_i(\Z)=\{f_i(x):\ x\in \Z\}\se\N$ for $i=1,2,3$, if
any $n\in\N$ can be written as $f_1(x)+f_2(y)+f_3(z)$ with $x,y,z\in \Z$ then we call the sum
$f_1(x)+f_2(y)+f_3(z)$ {\it universal over $\Z$}. For example, $p_3(x)+p_3(y)+p_3(z)$ is universal over $\Z$ by Gauss' result.

In 1862 Liouville (cf. \cite[p.\,82]{B}) determined all universal sums $ap_3(x)+bp_3(y)+cp_3(z)$ over $\Z$ with $a,b,c\in\Z^+=\{1,2,3,\ldots\}$.
Z.-W. Sun \cite{S07,S15} studied universal sums of the form $ap_i(x)+bp_j(y)+cp_k(z)$ with $a,b,c\in\N$ and $i,j,k\in\{3,4,\ldots\}$, where
$p_m(x)$ denotes the generalized polygonal number
$$(m-2)\bi x2+x=\f{x((m-2)x-(m-4))}2;$$ see also \cite{GPS,OS,GS,Oh,JOS} for subsequent work on some of Sun's conjectures posed in \cite{S07,S15}.
In 2017 Sun \cite{S17} investigated universal sums $x(ax+1)+y(by+1)+z(cz+1)$ over $\Z$ with $a,b,c\in\Z^+$,
and universal sums $x(ax+b)+y(ay+c)+z(az+d)$ over $\Z$ with $a,b,c,d\in\N$ and $a\gs b\gs c\gs d$. Quite recently, Sun \cite{S17a} considered more general questions
and investigated for what tuples $(a,b,c,d,e,f)$ with $a\gs c\gs e\gs 1$, $b\eq a\pmod 2$ and $0\ls b\ls a$, $d\eq c\pmod2$ and $0\ls d\ls c$,
$f\eq e\pmod2$ and $0\ls f\ls e$, the sum
$$\f{x(ax+b)}2+\f{y(cy+d)}2+\f{z(ez+f)}2$$
is universal over $\Z$. Such (ordered) tuples $(a,b,c,d,e,f)$ are said to be universal over $\Z$.
He showed such tuples with $b<a$, $d<c$, $f<e$, and $b\gs d$ if $a=c$, and $d\gs f$ if $c=e$,
must be in his list of 12082 candidates (cf. \cite[A286944]{S17b} and \cite{S18}), and conjectured that all such candidates are indeed universal over $\Z$.
Note that $$\{p_3(x):\ x\in\Z\}=\l\{\f{x(4x+2)}2=x(2x+1):\ x\in\Z\r\}.$$
Sun \cite{S17a} proved that some candidates $(a,b,c,d,e,f)$ are universal over $\Z$, e.g., $(5,1,3,1,1,1)$
(equivalent to $(5,1,4,2,3,1)$)
is universal over $\Z$. Sun even conjectured that any $n\in\N$ can be written as $x(x+1)/2+y(3y+1)/2+z(5z+1)/2$ with $x,y,z\in\N$.

In this paper, via the theory of ternary quadratic forms, we establish the universality (over $\Z$) of 44 concrete tuples $(a,b,c,d,e,f)$
on Sun's list of candidates.

\begin{theorem} \label{Th1.1} The tuples
\begin{align*}&(5,1,2,2,1,1),\ (6,0,3,3,3,1),\ (6,2,5,5,1,1),\ (6,6,3,3,3,1),
\\&  (8,2,3,1,1,1),\ (8,6,3,1,1,1),\ (8,8,3,1,1,1)
\end{align*}
are universal over $\Z$.
\end{theorem}
\begin{remark}\label{Rem1.1}
Sun \cite{S15} conjectured that any $n\in\N$ can be written as $p_3(x)+2p_3(y)+p_7(z)$ with $x,y,z\in\N$,
and J. Ju, B.-K. Oh and B. Seo \cite{JOS} proved that $p_3(x)+2p_3(y)+p_7(z)$ (or the tuple $(5,3,2,2,1,1)$)
is universal over $\Z$.
\end{remark}

\begin{theorem} \label{Th1.2} The tuples
\begin{align*}&(6,0,5,1,3,1),\ (6,0,5,3,3,1),\ (7,1,1,1,1,1),\ (7,1,2,0,1,1),
\\&(7,1,2,2,1,1),\ (7,1,3,1,1,1),\ (7,3,1,1,1,1),\ (7,3,2,0,1,1),
\\&(7,3,2,2,1,1),\ (7,3,3,1,1,1),\ (7,5,1,1,1,1),\ (7,5,3,1,1,1),
\\& (15,3,3,1,1,1),\ (15,5,1,1,1,1),\ (15,5,3,1,2,0),\ (15,5,3,1,2,2),
\\& (15,9,3,1,1,1),\ (21,7,3,1,2,2)
\end{align*}
are universal over $\Z$.
\end{theorem}
\begin{remark}\label{Rem1.2}
Our proof of Theorem \ref{Th1.2} involves the theory of genera of ternary quadratic forms.
Sun \cite{S15} conjectured that any $n\in\N$ can be written as $p_3(x)+y^2+p_9(z)$ (or $p_3(x)+2p_3(y)+p_9(z)$)
with $x,y,z\in\N$; along this line Ju, Oh and Seo \cite{JOS} proved that $p_3(x)+y^2+p_9(z)$
and $p_3(x)+2p_3(y)+p_9(z)$ are universal over $\Z$, i.e., the tuples $(7,5,2,0,1,1)$ and $(7,5,2,2,1,1)$ are universal over $\Z$.\end{remark}

\begin{theorem} \label{Th1.3} {\rm (i)} The tuples $(5,5,3,1,3,1)$, $(5,5,3,3,3,1)$, $(6,4,5,5,1,1)$ and $(7,7,3,1,1,1)$ are universal over $\Z$.

{\rm (ii)} All the five tuples
$$(6,2,5,1,1,1),\ (6,2,5,5,1,1),\ (6,4,5,1,1,1),\ (15,5,6,2,1,1),\ (15,5,6,4,1,1)$$
are universal over $\Z$.
\end{theorem}
\begin{remark}\label{Rem1.3}
Our proof of Theorem \ref{Th1.3}(i) employs the Minkowski-Siegel formula (cf. \cite [pp.\, 173-174]{Ki}).
Sun \cite{S15} conjectured that any $n\in\N$ can be written as $p_3(x)+p_7(y)+2p_5(z)$ (or $p_3(x)+p_7(y)+p_8(z)$)
with $x,y,z\in\N$; along this line Ju, Oh and Seo \cite{JOS} proved that $p_3(x)+p_7(y)+2p_5(z)$
and $p_3(x)+p_7(y)+p_8(z)$ are universal over $\Z$, i.e., the tuples $(6,2,5,3,1,1)$ and $(6,4,5,3,1,1)$ are universal over $\Z$.\end{remark}

Similarly to \cite[Theorem 1.4]{S17a}, we observe that
\begin{equation}\label{1.1}\{p_3(x)+p_5(y):\ x,y\in\Z\}=\{p_5(x)+3p_5(y):\ x,y\in\Z\}.
\end{equation}
In fact,
\begin{align*}&n\in\{p_3(x)+p_5(y):\ x,y\in\Z\}
\\\iff& 24n+4\in\{3(2x+1)^2+(6y-1)^2:\ x,y\in\Z\}
\\\iff& 24n+4\in\{3u^2+v^2:\ u,v\in\Z\ \&\ 2\nmid uv\}
\end{align*}
and
\begin{align*}&n\in\{3p_5(x)+p_5(y):\ x,y\in\Z\}
\\\iff& 24n+4\in\{3(6x-1)^2+(6y-1)^2:\ x,y\in\Z\}
\\\iff& 24n+4\in\{3u^2+v^2:\ u,v\in\Z,\ 2\nmid uv\ \&\ 3\nmid u\}.
\end{align*}
If $u$ and $v$ are odd integers with $3\mid u$ and $3\nmid v$, then
$$3u^2+v^2=3\l(\f{u\pm v}2\r)^2+\l(\f{3u\mp v}2\r)^2$$
with $(u\pm v)/2$ not divisible by $3$. Therefore (\ref{1.1}) holds.
In view of (1.1) and Theorems 1.1-1.3, we have the following consequence.

\begin{corollary}\label{Cor1.1} The tuples
\begin{align*}&(9,3,7,1,3,1),\ (9,3,7,3,3,1),\ (9,3,7,5,3,1),
\\& (9,3,7,7,3,1),\ (9,3,8,2,3,1),\ (9,3,8,6,3,1),
\\& (9,3,8,8,3,1),\ (15,3,9,3,3,1),\ (15,9,9,3,3,1)
\end{align*}
are universal over $\Z$.
\end{corollary}

\begin{theorem}\label{Th1.4} The tuple $(16,4,2,0,1,1)$ is universal over $\Z$. In other words,  any $n\in\N$
can be written as $p_3(x)+y^2+2z(4z+1)$ with $x,y,z\in\Z$.
\end{theorem}
\begin{remark}\label{Rem1.4} This result is closely related to the form $x^2+y^2+32z^2$.
Sun \cite{S17a} even conjectured that any $n\in\N$ can be written as $p_3(x)+y^2+2z(4z-1)$ with $x,y,z\in\N$.
\end{remark}

We will show Theorems 1.1-1.4 in Sections 2-5 respectively.

\section{Proof of Theorem \ref{Th1.1}}

\begin{lemma}\label{Lem2.2} {\rm (i)} For any $n\in\N$, we can write $12n+5$ as $x^2+y^2+(6z)^2$ with $x,y,z\in\Z$.

{\rm (ii)} Let $n\in\Z^+$ and $\da\in\{0,1\}$. Then we can write $6n+1$ as $x^2+3y^2+6z^2$ with $x,y,z\in\Z$ and $x\eq\da\pmod2$.
\end{lemma}

\begin{remark}\label{Rem2.2} Lemma \ref{Lem2.2} is a known result due to the second author, see \cite[Theorem 1.7(iii) and Lemma 3.3]{S15} and \cite[Remark 3.1]{S17}.
\end{remark}

The following lemma (cf. \cite [pp.\,12-14]{Jagy14}) occurred in a 1993 letter of J.S. Hsia to I. Kaplansky.

\begin{lemma}\label{Lem2.3} For each $n \in \N$, we can write $6n+5$ as $x^2+y^2+10z^2$ with $x,y,z\in \Z$.
\end{lemma}

For $a,b,c\in\Z^+$, we define
$$E(a,b,c)=\{n\in\N:\ n\not=ax^2+by^2+cz^2\ \t{for any}\ x,y,z\in\Z\}.$$
L.E. Dickson \cite[pp.\,112-113]{D39} listed 102 diagonal quadratic forms $ax^2+by^2+cz^2$
for which the structure of $E(a,b,c)$ is known explicitly. For example, the Gauss-Legendre theorem asserts that
$E(1,1,1)=\{4^k(8l+7):\ k,l\in\N\}$.

In 1996 W. Jagy \cite{Jagy96} showed the following result (cf. \cite [pp.\,25-26]{Jagy14}).
\begin{lemma}\label{Lem2.4} We have
$$E(1,4,9)=\{2\}\cup\bigcup_{k,l\in\N}\{4^k(8l+7),\ 8l+3,\ 9l+3\}.$$
\end{lemma}

\medskip
\noindent{\it Proof of Theorem \ref{Th1.1}.} (i) We want to prove the universality of
$(5,r,2,2,1,1)$ over $\Z$ for $r\in\{1,3\}$. Let $n\in\N$. Clearly,
\begin{align*}&n=p_3(x)+y(y+1)+\f{z(5z+r)}2
\\\iff&40n+r^2+15=5(2x+1)^2+10(2y+1)^2+(10z+r)^2.
\end{align*}
Since
$$E(1,5,10)=\{25^lm:\ l,m\in\N\ \mbox{and}\ m\eq 2,3\pmod5\}$$
by Dickson \cite[pp.\,112-113]{D39},
we have $40n+r^2+15\in\{u^2+5v^2+10w^2:\ u,v,w\in\N\}$. Thus we can write
$$40n+r^2+15=(2^kx_0)^2+5(2^ky_0)^2+10(2^kz_0)^2=4^k(x_0^2+5y_0^2+10z_0^2)$$
with $k\in\N$, $x_0,y_0,z_0\in\Z$, and $x_0,y_0,z_0$ not all even. In the case $k=0$, if $2\mid z_0$ then $x_0^2+5y_0^2\eq r^2+15\eq0\pmod8$ and hence
$x_0\eq y_0\eq0\pmod2$ which contradicts that $x_0,y_0,z_0$ are not all even, thus $2\nmid z_0$ and also $2\nmid x_0y_0$ since
$x_0^2+5y_0^2\eq r^2+15-10z_0^2\eq6\pmod8$.

 It is easy to verify the following new identity:
 \begin{equation}\label{2.1}
 4^2(x^2+5y^2+10z^2)=(x\pm 5y-10z)^2+5(x\mp 3y-2z)^2+10(x\pm y+2z)^2.
 \end{equation}
 If $x,y,z$ are odd integers, then $(x+\ve y)/2+z$ is odd for some $\ve\in\{\pm1\}$, hence by (\ref{2.1}) we have
 $$4(x^2+5y^2+10z^2)=\tilde x^2+5\tilde y^2+10\tilde z^2$$
 with
 $$\tilde x=\f{x+\ve y}2+2\ve y-5z,\ \tilde y=\f{x+\ve y}2-2\ve y-z,\ \tilde z=\f{x+\ve y}2+z$$
 all odd. Thus, if $2\nmid x_0y_0z_0$ then
 \begin{equation}\label{2.2}40n+r^2+15=4^k(x_0^2+5y_0^2+10z_0^2)\in\{x^2+5y^2+10z^2:\ x_,y,z\ \mbox{are odd}\}.
 \end{equation}

 If $x_0\not\eq y_0\pmod2$, then $x_0^2+5y_0^2+10z_0^2\eq1\pmod2$ and $k\gs2$ since $40n+r^2+15\eq0\pmod8$, hence by (\ref{2.1}) we have
 $$4^2(x_0^2+5y_0^2+10z_0^2)=\bar x_0^2+5\bar y_0^2+10\bar z_0^2$$
 with $\bar x_0=x_0-5y_0-10z_0$, $\bar y_0=x_0+3y_0-2z_0$ and $\bar z_0=x_0-y_0+2z_0$ all odd, and therefore (\ref{2.2}) holds.

 Now we suppose that $k>0$, $2\mid x_0y_0z_0$ and $x_0\eq y_0\pmod2$. By (\ref{2.1}),
 $$4(x_0^2+5y_0^2+10z_0^2)=x_1^2+5y_1^2+10z_1^2$$
 with
 $$x_1=\f{x_0-y_0}2-2y_0-5z_0,\ y_1=\f{x_0-y_0}2+2y_0-z_0,\ z_1=\f{x_0-y_0}2+z_0$$
 If $x_0$ and $y_0$ are odd, then we may assume $x_0\not\eq y_0-2z_0\pmod4$ without loss of generality (otherwise we replace $x_0$ by $-x_0$), and hence
 $x_1,y_1,z_1$ are all odd. If $x_0,y_0,(x_0-y_0)/2$ are all even, then $z_0$ is odd and so are $x_1,y_1,z_1$.
 If $x_0$ and $y_0$ are even with $x_0\not\eq y_0\pmod4$, then $z_0$ is odd and we may assume $z_0\eq (y_0-x_0)/2\pmod4$ without loss of generality (otherwise we replace $z_0$ by $-z_0$), hence $z_1\eq0\pmod4$, $y_1=z_1+2(y_0-z_0)\eq0\pmod2$ and $(x_1-y_1)/4\eq -y_0-z_0\eq1\pmod2$, therefore by (\ref{2.1}) we have
 $$x_1^2+5y_1^2+10z_1^2=x_2^2+5y_2^2+10z_2^2$$
 with
 $$x_2=\f{x_1-5y_1-10z_1}4,\ y_2=\f{x_1+3y_1-2z_1}4,\ z_2=\f{x_1-y_1+2z_1}4$$
 all odd. So we still have (\ref{2.2}).

 By the above, there always exist odd integers $x,y,z$ such that $40n+r^2+15=x^2+5y^2+10z^2$. Write $y=2u+1$ and $z=2v+1$ with $u,v\in\Z$.
 As $x^2\eq r^2\pmod5$, either $x$ or $-x$ has the form $10w+r$ with $w\in\Z$. Therefore
 $$40n+r^2+15=(10w+r)^2+5(2u+1)^2+10(2v+1)^2$$
 and hence $n=p_3(u)+v(v+1)+w(5w+r)/2$. This proves the universality of $(5,r,2,2,1,1)$ over $\Z$.

(ii) Let $n\in\N$ and $r\in\{1,3\}$. It is easy to see that
\begin{align*}&n=p_3(x)+\f{y(3y+1)}2+z(4z+r)
\\\iff&48n+3r^2+8=6(2x+1)^2+2(6y+1)^2+3(8z+r)^2.
\end{align*}

Since
$$E(2,3,6)=\{3q+1:\ q\in\N\}\cup\{4^k(8l+7):\ k,l\in\N\}$$
by Dickson \cite[pp.\,112-113]{D39}, we see that $48n+3r^2+8=2x^2+3y^2+6z^2$ for some $x,y,z\in\Z$.
Clearly, $y^2+2z^2\not=0$, and hence by \cite[Lemma 2.1]{S15} we have $y^2+2z^2=y_0^2+2z_0^2$ for some $y_0,z_0\in\Z$ not all divisible by 3.
Thus, without any loss of generality, we simply assume that $3\nmid y$ or $3\nmid z$. Note that $3\nmid x$, $2\nmid y$, and $x\eq z\pmod 2$
since $2(x^2+z^2)\eq 2x^2+6z^2\eq 3r^2+8-3y^2\eq0\pmod4$. If $3\mid y$ and $3\nmid z$, then $z$ or $-z$ is congruent to $x+y$ modulo 3.
If $3\nmid y$ and $3\mid z$, then $y$ or $-y$ is congruent to $x+z$ modulo 3. If $3\nmid yz$, then $\ve_1 y\eq \ve_2 z\eq x\pmod3$
for some $\ve_1,\ve_2\in\{\pm1\}$. So, without loss of generality, we may assume that $x+y+z\eq0\pmod3$
(otherwise we may change signs of $x,y,z$ suitably). Note that
$$48n+3r^2+8=2x^2+3y^2+6z^2=2a^2+3b^2+6c^2,$$
where $a=y+z$, $b=(2x-y+2z)/3$ and $c=(x+y-2z)/3$ are integers.
If $x\eq z\eq1\pmod2$, then $x,y,z$ are all odd. If $x\eq z\eq0\pmod2$, then $a,b,c$ are all odd.

By the above, $48n+3r^2+8=2a^2+3b^2+6c^2$ for some odd integers $a,b,c$.
Since $3b^2\eq 3r^2+8-2a^2-6c^2\eq3r^2\pmod{16}$, we can write $b$ or $-b$ as $8w+r$ with $w\in\Z$.
Clearly, $a$ or $-a$ has the form $6u+1$ with $u\in\Z$, and $c=2v+1$ for some $v\in\Z$.
Therefore
$$48n+3r^2+8=2(6u+1)^2+3(8w+r)^2+6(2v+1)^2$$
and hence
$n=u(3u+1)/2+p_3(v)+w(4w+r)$. This proves the universality of $(8,2r,3,1,1,1)$ over $\Z$.

(iii) Let $n\in\N$.  By Lemma \ref{Lem2.2}(ii), we can write $6n+7$ in the form $x^2+3y^2+6z^2$ with $x,y,z\in\Z$ and  $x\equiv n+1\pmod 2$. Clearly, $y\eq n\pmod2$.
Since $6z^2\eq 6n+7-(n+1)^2-3n^2\equiv 6\pmod 4$, we have $2\nmid z$.
Hence
\begin{equation*}
24n+28=4(6n+7)=4(x^2+3y^2+6z^2)=(x-3y)^2+3(x+y)^2+24z^2
\end{equation*}
with $x-2y$, $x+2y$ and $z$ all odd. Note that $x-3y$ or $3y-x$ has the form $6w+1$ with $w\in\Z$.
Write $x+y=2u+1$ and $z=2v+1$ with $u,v\in\Z$. Then
$$24n+28=(6w+1)^2+3(2u+1)^2+24(2v+1)^2$$
and hence $n=w(3w+1)/2+p_3(u)+8p_3(v)$. This proves the universality of $(8,8,3,1,1,1)$.

(iv) Let $n\in\N$. By Lemma \ref{Lem2.3}, we can write $6n+5$ as $x^2+y^2+10z^2$ with $x,y,z\in\Z$. Clearly, $x\not\eq y\pmod2$.
Since $x^2+y^2+z^2\eq2\pmod 3$, exactly one of $x,y,z$ is divisible by $3$.
Without loss of generality, we may assume that $x+y+z\eq0\pmod3$
(otherwise we adjust signs of $x,y,z$ suitably to meet our purpose). Observe that
$$4(x^2+y^2+10z^2)=2(x-y)^2+3\l(\f{x+y+10z}3\r)^2+15\l(\f{x+y-2z}3\r)^2.$$
So, $4(6n+5)=2a^2+3b^2+15c^2$ for some odd integers $a,b,c$. As $3\nmid a$, we may write $a$ or $-a$ as $6w+1$ with $w\in\Z$.
Write $b=2u+1$ and $c=2v+1$ with $u,v\in\Z$. Then
$$24n+20=2(6w+1)^2+3(2u+1)^2+15(2v+1)^2$$
and so $n=p_3(u)+5p_3(v)+w(3w+1)$. This proves the universality of $(6,2,5,5,1,1)$ over $\Z$.

(v) Let $n\in\N$. By Lemma \ref{Lem2.2}(i), we can write $12n+5$ in the form $x^2+y^2+(6z)^2$ with $x,y,z\in\Z$.
It follows that $24n+10=(x+y)^2+(x-y)^2+72z^2$. As $(x+y)^2+(x-y)^2\eq10\eq2\pmod4$, both $x+y$ and $x-y$ are odd.
Since $(x+y)^2+(x-y)^2\eq 10\eq1\pmod3$, exactly one of $x+y$ and $x-y$ is divisible by 3. So $(x+y)^2+(x-y)^2=(6u+1)^2+(6v+3)^2$
for some $u,v\in\Z$. Therefore
$$24n+10=(6u+1)^2+(6v+3)^2+72z^2,$$
i.e., $n=u(3u+1)/2+3p_3(v)+3z^2$. This proves the universality of $(6,0,3,3,3,1)$ over $\Z$.

By Lemma \ref{Lem2.4}, we can write $12n+14$ in the form $x^2+4y^2+9z^2$ with $x,y,z\in\Z$.
Since $x^2+z^2\eq 14\pmod4$, we have $2\nmid xz$. Observe that
\begin{equation*}
24n+28=2(x^2+4y^2+9z^2)=(x-2y)^2+(x+2y)^2+18z^2
\end{equation*}
with $x\pm 2y$ and $z$ all odd.
Clearly, exactly one of $x-2y$ and $x+2y$ is divisible by 3.
So, for some $u,v,w\in\Z$ we have
$$24n+28=(6u+1)^2+9(2v+1)^2+18(2w+1)^2$$
and hence $n=u(3u+1)/2+3p_3(v)+6p_3(w)$. This proves the universality of $(6,6,3,3,3,1)$ over $\Z$.

The proof of Theorem 1.1 is now complete.
\qed

\section{Proof of Theorem \ref{Th1.2}}

The following lemma is one of the most important theorems about integral representations of quadratic forms (cf. \cite [p.\,129]{C}).
\begin{lemma}\label{Lem3.1} Let $f$ be a nonsingular integral quadratic form and let $m$ be a nonzero integer which is represented by $f$ over the real field $\R$
and the ring $\Z_p$ of $p$-adic integers for each prime $p$.
Then $m$ is represented by some form $f^*$ over $\Z$ where $f^*$ is in the genus of $f$.
\end{lemma}

\begin{lemma}\label{Lem3.2} {\rm (i)} \cite[Lemma 3.2]{S15} If $x^2+3y^2\eq4\pmod 8$ with $x,y\in\Z$, then $x^2+3y^2=u^2+3v^2$
for some odd integers $u$ and $v$.

{\rm (ii)}  \cite[Lemma 3.6]{S15} If $w=x^2+7y^2>0$ with $x,y\in\Z$ and $8\mid w$,
then $w=u^2+7v^2$ for some odd integers $u$ and $v$.

{\rm (iii)} \cite[Lemma 5.1]{S17a} If $w=3x^2+5y^2>0$ with $x,y\in\Z$ and $8\mid w$,
then $w=3u^2+5v^2$ for some odd integers $u$ and $v$.
\end{lemma}

\medskip
\noindent{\it Proof of Theorem \ref{Th1.2}}. (i) Let $n\in\N$. Clearly,
$$n=p_3(x)+p_3(y)+5z(3z+1)/2\iff 24n+11=3(2x+1)^2+3(2y+1)^2+5(6z+1)^2.$$

There are two classes in the genus of $3x^2+3y^2+5z^2$, and the one not containing $3x^2+3y^2+5z^2$ has the representative
\begin{align*}3x^2+2y^2+8z^2-2yz=&3x^2+3\l(\f y2+z\r)^2+5\l(\f y2-z\r)^2
\\=&3x^2+3\l(\f{y-3z}2\r)^2+5\l(\f{y+z}2\r)^2
\end{align*}
If $24n+11=3x^2+2y^2+8z^2-2yz$ with $y$ odd and $z$ even, then $3x^2\eq 11-2y^2\eq9\pmod 4$ which is impossible.
Thus, if $24n+11\in\{3x^2+2y^2+8z^2-2yz:\ x,y,z\in\Z\}$ then $24n+11\in\{3x^2+3y^2+5z^2:\ x,y,z\in\Z\}$.
With the help of Lemma \ref{Lem3.1}, there are $x,y,z\in\Z$ such that $24n+11=3x^2+3y^2+5z^2$.
As $5z^2\not\eq11\pmod4$, $x$ and $y$ cannot be both even. Without loss of generality, we assume that $2\nmid x$.
Then $3y^2+5z^2\eq 11-3x^2\eq0\pmod 8$ and $3y^2+5z^2\not=0$. By Lemma \ref{Lem3.2}(iii), $3y^2+5z^2=3y_0^2+5z_0^2$ for some odd integers $y_0$ and $z_0$.
Write $x=2u+1$ and $y_0=2v+1$ with $u,v\in\Z$. As $2\nmid z_0$ and $3\nmid z_0$, $z_0$ or $-z_0$ has the form $6w+1$ with $w\in\Z$. Thus
$24n+11=3(2u+1)^2+3(2v+1)^2+5(6w+1)^2$ and hence $n=p_3(u)+p_3(v)+5w(3w+1)/2$. This proves the universality of $(15,5,1,1,1,1)$ over $\Z$.
\medskip

(ii) Let $n\in\N$ and $r\in\{1,3\}$. Obviously,
\begin{align*}&n=p_3(x)+\f{y(3y+1)}2+3\f{z(5z+r)}2
\\\iff&120n+9r^2+20=15(2x+1)^2+5(6y+1)^2+9(10z+r)^2.
\end{align*}

There are two classes in the genus of $x^2+15y^2+5z^2$, and the one not containing $x^2+15y^2+5z^2$ has the representative
\begin{align*}4x^2+4y^2+5z^2+2xy=&\l(\f x2+2y\r)^2+15\l(\f x2\r)^2+5z^2
\\=&\l(2x+\f y2\r)^2+15\l(\f y2\r)^2+5z^2.
\end{align*}
If $120n+9r^2+20=4x^2+4y^2+5z^2+2xy$ with $x,y,z\in\Z$, then $2xy\eq 9r^2-5z^2\eq0\pmod4$ and hence $x$ or $y$ is even.
Thus, with the help of Lemma \ref{Lem3.1}, we can always write
$120n+9r^2+20=x^2+15y^2+5z^2$ with $x,y,z\in\Z$. Since $x^2+5z^2\eq 20\eq2\pmod3$, $x=3x_0$ for some $x_0\in\Z$.
As $15y^2\not\eq 9r^2\pmod4$, $x$ and $z$ cannot be both even.
If $2\nmid x$, then $5(3y^2+z^2)\eq 9r^2+20-x^2\eq4\pmod 8$ and hence by Lemma \ref{Lem3.2}(i) we can write $3y^2+z^2$ as $3y_0^2+z_0^2$ with $y_0$ and $z_0$ both odd.
If $2\nmid z$, then $x^2+15y^2\not=0$ and $x^2+15y^2=3(3x_0^2+5y^2)\eq 9r^2+20-5z^2\eq0\pmod8$, hence by Lemma \ref{Lem3.2}(iii) we can write
$3x_0^2+5y^2$ as $3x_1^2+5y_1^2$ with $x_1$ and $y_1$ both odd.

By the above, there are odd integers $x,y,z$ such that $120n+9r^2+20=9x^2+15y^2+5z^2$.
Write $y=2u+1$ with $u\in\Z$.
As $3\nmid z$, we can write $z$ or $-z$ as $6v+1$ with $v\in\Z$. Since $x^2\eq r^2\pmod 5$, we can write
$x$ or $-x$ as $10w+r$ with $w\in\Z$. Thus
$$120n+9r^2+20=15(2u+1)^2+5(6v+1)^2+9(10z+r)^2$$
and hence
$n=p_3(x)+y(3y+1)/2+3z(5z+r)/2$ with $x,y,z\in\Z$.
This proves the universality of $(15,3r,3,1,1,1)$ over $\Z$.

(iii) Let $n\in\N$ and $r\in\{1,3\}$. Obviously,
\begin{align*}&n=3x^2+\f{y(3y+1)}2+\f{z(5z+r)}2
\\\iff&120n+3r^2+5=360x^2+5(6y+s)^2+3(10z+r)^2.
\end{align*}
If $60n+(3r^2+5)/2=4x^2+4y^2+5z^2+2xy$ with $x,y,z\in\Z$, then $x$ or $y$ must be even.
Thus, as in part (ii), $60n+(3r^2+5)/2=x^2+5y^2+15z^2$ for some $x,y,z\in\Z$.
Note that $x^2+y^2\eq z^2\pmod4$.  If $y$ is odd, then $2\mid x$, $2\nmid z$ and we may assume $y\not\eq z\pmod4$ (otherwise it suffices to change the sign of $z$),
hence
$$y^2+3z^2=\l(\f{y-3z}2\r)^2+3\l(\f{y+z}2\r)^2$$
with $y_1=(y-3z)/2$ and $z_1=(y+z)/2$ both even.
So, without loss of generality, we may simply assume that $2\mid y$ and $x\eq z\pmod2$.
Observe that
$$120n+3r^2+5=2(x^2+5y^2+15z^2)=3a^2+5b^2+10y^2.$$
with $a=(x+5z)/2$ and $b=(x-3z)/2$ both integral. Since $3a^2+5b^2\eq 5s^2+3t^2-10y^2\eq0\pmod8$ and $3a^2+5b^2>0$, by Lemma \ref{Lem3.2}(iii)
we can write $3a^2+5b^2=3c^2+5d^2$ with $c$ and $d$ both odd. Thus
$$120n+3r^2+5=3c^2+5d^2+40\l(\f y2\r)^2.$$
As $(y/2)^2\eq 5(1-d^2)\eq d^2-1\pmod3$, we must have $3\nmid d$ and $3\mid y$.
Write $y=6u$ with $u\in\Z$. Clearly, $d$ or $-d$ has the form $6v+1$ with $v\in\Z$. Since $c^2\eq r^2\pmod5$,
we may write $c$ or $-c$ as $10w+r$ with $w\in\Z$. Therefore
$$120n+3r^2+5=3(10w+r)^2+5(6v+1)^2+40(3u)^2$$
and hence $n=3u^2+v(3v+1)/2+w(5w+r)/2$. This proves the universality of $(6,0,5,r,3,1)$ over $\Z$.

(iv) Let $n\in\N$ and $\da\in\{0,1\}$. Clearly,
\begin{align*}&n=x(x+\da)+\f{y(3y+1)}2+5\f{z(3z+1)}2
\\\iff&24n+6(\da+1)=6(2x+\da)^2+(6y+1)^2+5(6z+1)^2.
\end{align*}

There are two classes in the genus of $x^2+5y^2+6z^2$, and the one not containing $x^2+5y^2+6z^2$ has the representative $3x^2+3y^2+4z^2-2yz+2zx$.
If $24n+6(\da+1)=3x^2+3y^2+4z^2-2yz+2zx$, then $u=(x+y)/2$ and $v=(x-y)/2$ are integers, and
$$24n+6(\da+1)=6u^2+6v^2+4z^2+4vz=6u^2+5v^2+(v+2z)^2.$$
Thus, by Lemma \ref{Lem3.1}, $24n+6(\da+1)=x^2+5y^2+6z^2$ for some $x,y,z\in\Z$.
Since $x^2\eq -5y^2\eq y^2\pmod3$, we may assume that $x\eq y\pmod3$ without loss of generality.
If $z\not\eq\da\pmod2$, then $x^2+5y^2\eq 6(\da+1)-6z^2\eq 6(\da+1)-6(1-\da)\eq4\da\pmod8$, hence both $x$ and $y$ are even and $(x-y)/2\eq\da\pmod2$,
and thus
$$x^2+5y^2+6z^2=\l(z-\f{5(x-y)}6\r)^2+5\l(\f{x-y}6+z\r)^2+6\l(\f{x-y}6+y\r)^2$$
with $(x-y)/6+y\eq (x-y)/2\eq\da\pmod2$.

By the above, $24n+6(\da+1)=x^2+5y^2+6z^2$ for some $x,y,z\in\Z$ with $x,y,z\in\Z$ with $z\eq\da\pmod2$.
Since $x^2+5y^2$ is a positive multiple of 3, by \cite[Lemma 2.1]{S15} we can write $x^2+5y^2=x_0^2+5y_0^2$ with $x_0y_0\in\Z$ and $3\nmid x_0y_0$.
So, there are $x,y,z\in\Z$ with $x\eq y\not\eq0\pmod3$ and $z\eq\da\pmod2$
such that $24n+6(\da+1)=x^2+5y^2+6z^2$. Write $z=2w+\da$ with $w\in\Z$.
Since $x^2+5y^2\eq6\pmod8$, both $x$ and $y$ are odd. Thus $x$ or $-x$ has the form $6u+1$ with $u\in\Z$,
and $y$ or $-y$ has the form $6v+1$ with $v\in\Z$. Therefore
$$24n+6(\da+1)=(6u+1)^2+5(6v+1)^2+6(2w+\da)^2$$
and hence $n=w(w+\da)+u(3u+1)/2+5v(3v+1)/2$.
This proves the universality of $(15,5,3,1,2,2\da)$ over $\Z$.

(v) Let $n\in\N$. Apparently,
\begin{align*}&n=x(x+1)+\f{y(3y+1)}2+7\f{z(3z+1)}2
\\\iff&24n+14=6(2x+1)^2+(6y+1)^2+7(6z+1)^2.
\end{align*}

There are two classes in the genus of $x^2+6y^2+7z^2$, and the one not containing $x^2+6y^2+7z^2$ has the representative
$$2x^2+5y^2+5z^2-4yz=2x^2+10u^2+10v^2-4(u+v)(u-v)=2x^2+6u^2+14v^2$$
with $u=(y+z)/2$ and $v=(y-z)/2$.
If $24n+14=2x^2+6u^2+14v^2$ for some $x,u,v\in\Z$ with $x\not\eq v\pmod2$,
then $14\eq 2+6u^2\pmod 8$ which is impossible. If $24n+14=2x^2+6u^2+14v^2$ with $x,u,v\in\Z$ and $x\eq v\pmod2$,
then
$$24n+14=6u^2+\l(\f{x-7v}2\r)^2+7\l(\f{x+v}2\r)^2.$$

By the above and Lemma \ref{Lem3.1}, there are $x,y,z\in\Z$ such that $24n+14=6x^2+y^2+7z^2$.
If $2\mid x$, then $y^2+7z^2\eq 6-6x^2\eq6\pmod8$ which is impossible. So $x=2u+1$ for some $u\in\Z$.
Note that $y^2+7z^2\eq 6-6x^2\eq0\pmod8$ and $y^2+7z^2\not=0$.
 Applying Lemma \ref{Lem3.2}(ii)
we can write $y^2+7z^2$ as $y_0^2+7z_0^2$ with $y_0$ and $z_0$ both odd.
Note that $y_0^2+z_0^2\eq y_0^2+7z_0^2\eq 14\eq2\pmod 3$.
So $y_0$ or $-y_0$ can be written as $6v+1$ with $v\in\Z$, and $z_0$ or $-z_0$ has the form $6w+1$ with $w\in\Z$.
Thus
$$24n+14=6x^2+y_0^2+7z_0^2=6(2u+1)^2+(6v+1)^2+7(6w+1)^2$$
and hence $n=u(u+1)+v(3v+1)/2+7z(3z+1)/2$. This proves the universality of $(21,7,3,1,2,2)$.

(vi) Let $r\in\{1,3,5\}$ and $n\in\N$. Clearly,
$$n=p_3(x)+p_3(y)+\f{z(7z+r)}2\iff 56n+14+r^2=7(2x+1)^2+7(2y+1)^2+(14z+r)^2.$$

There are two classes in the genus of $x^2+7y^2+7z^2$, and the one not containing $x^2+7y^2+7z^2$ has the representative
\begin{align*}2x^2+4y^2+7z^2+2xy=&\l(\f x2+2y\r)^2+7\l(\f x2\r)^2+7z^2.
\\=&\l(\f{x-3y}2\r)^2+7\l(\f{x+y}2\r)^2+7z^2
\end{align*}
If $56n+14+r^2=2x^2+4y^2+7z^2+2xy$ with $x$ odd and $y$ even, then $15\eq 14+r^2\eq 2x^2+7z^2\eq 9\pmod 4$ which is impossible.
Thus, if $56n+14+r^2\in\{2x^2+4y^2+7z^2+2xy:\ x,y,z\in\Z\}$ then $56n+14+r^2\in\{x^2+7y^2+7z^2:\ x,y,z\in\Z\}$.
With the help of Lemma \ref{Lem3.1}, there are $x,y,z\in\Z$ such that $56n+14+r^2=x^2+7y^2+7z^2$.
As $x^2\not\eq 14+r^2\eq15\pmod4$, $y$ and $z$ cannot be both even. Without loss of generality, we assume that $2\nmid z$.
Then $x^2+7y^2\eq 14+r^2-7z^2\eq0\pmod 8$ and $x^2+7y^2\not=0$. By Lemma \ref{Lem3.2}(ii), $x^2+7y^2=x_0^2+7y_0^2$ for some odd integers $x_0$ and $y_0$.
Now $56n+14+r^2=x_0^2+7y_0^2+7z^2$. Clearly, $x_0$ or $-x_0$ has the form $14w+r$ with $w\in\Z$.
Write $y_0=2u+1$ and $z=2v+1$ with $u,v\in\Z$. Then
$$56n+14+r^2=(14w+r)^2+7(2u+1)^2+7(2v+1)^2$$
and hence $n=p_3(u)+p_3(v)+w(7w+r)/2$. This proves the universality of $(7,r,1,1,1,1)$ over $\Z$.

(vii) Let $n\in\N$ and $t\in\{1,3,5\}$. Clearly,
\begin{align*}&n=p_3(x)+\f{y(3y+1)}2+\f{z(7z+t)}2
\\\iff&168n+28+3t^2=21(2x+1)^2+7(6y+1)^2+3(14z+t)^2.
\end{align*}

There are two classes in the genus of $3x^2+21y^2+7z^2$, and the one not containing $3x^2+21y^2+7z^2$ has the representative
\begin{align*}6x^2+12y^2+7z^2+6xy=&3\l(\f x2+2y\r)^2+21\l(\f x2\r)^2+7z^2.
\\=&3\l(\f{x-3y}2\r)^2+21\l(\f{x+y}2\r)^2+7z^2
\end{align*}
If $168n+28+3t^2=6x^2+12y^2+7z^2+6xy$ with $x$ odd and $y$ even, then $31\eq 28+3t^2\eq 6x^2+7z^2\eq 13\pmod 4$ which is impossible.
Thus, if $168n+28+3t^2\in\{6x^2+12y^2+7z^2+6xy:\ x,y,z\in\Z\}$ then $168n+28+3t^2\in\{3x^2+21y^2+7z^2:\ x,y,z\in\Z\}$.
With the help of Lemma \ref{Lem3.1}, there are $x,y,z\in\Z$ such that $168n+28+3t^2=3x^2+21y^2+7z^2$.
As $21y^2\not\eq 28+3t^2\eq31\pmod4$, $x$ and $z$ cannot be both even.
If $2\nmid x$, then $21y^2+7z^2\eq28+3t^2-3x^2\eq4\pmod 8$ and hence by Lemma \ref{Lem3.2}(i) we can write $3y^2+z^2$ as $3y_0^2+z_0^2$ with $y_0,z_0$ odd integers.
Note that $x^2+7y^2\not=0$ since $7\nmid t$.
If $2\nmid z$, then $3(x^2+7y^2)\eq 28+3t^2-7z^2\eq0\pmod 8$ and hence by Lemma \ref{Lem3.2}(ii) $x^2+7y^2=x_0^2+7y_0^2$ for some odd integers $x_0$ and $y_0$.

By the above, there are odd integers $x,y,z$ such that $168n+28+3t^2=3x^2+7y^2+21z^2$.
Write $z=2u+1$ with $u\in\Z$.
As $y^2\eq 1\pmod3$, $y$ or $-y$ has the form $6v+1$ with $v\in\Z$.
Since $x^2\eq t^2\pmod 7$, $x$ or $-x$ has the form $14w+t$ with $w\in\Z$.
Thus
$$168n+28+3t^2=3(14w+t)^2+7(6v+1)^2+21(2u+1)^2$$
and hence
$n=p_3(u)+v(3v+1)/2+w(7w+t)/2$. This proves the universality of $(7,t,3,1,1,1)$ over $\Z$.

(viii) Let $\da\in\{0,1\}$ and $r\in\{1,3,5\}$.  Clearly,
\begin{align*}&n=p_3(x)+y(y+\da)+\f{z(7z+r)}2
\\\iff&56n+14\da+r^2+7=7(2x+1)^2+14(2y+\da)^2+(14z+r)^2.
\end{align*}

There are two classes in the genus of $x^2+7y^2+14z^2$, the one not containing $x^2+7y^2+14z^2$ has the representative
$$2x^2+7y^2+7z^2=2x^2+14\l(\f{y+z}2\r)^2+14\l(\f{y-z}2\r)^2.$$
If $56n+14\da+r^2+7=2x^2+14y^2+14z^2$ with $x,y,z\in\Z$ and $y,z\not\eq x\pmod 2$, then
$2x^2\eq14\da+r^2+7\eq 2\da\pmod 4$, hence $x^2\eq\da\pmod 4$ and also $y\eq z\eq\da\pmod2$
since
$$-2(y^2+z^2)\eq 14(y^2+z^2)\eq 14\da+r^2+7-2\da\eq-4\da\pmod8,$$
this contradicts with $y,z\not\eq x\pmod 2$.
If $56n+14\da+r^2+7=2x^2+14y^2+14z^2$ with $x,y,z\in\Z$ and $x\eq y\pmod2$, then
$$56n+14\da+r^2+7=\l(\f{x-7y}2\r)^2+7\l(\f{x+y}2\r)^2+14z^2.$$

In view of Lemma \ref{Lem3.1} and the above, there are $x,y,z\in\Z$ such that $56n+14\da+r^2+7=x^2+7y^2+14z^2$.
If $z\not\eq\da\pmod2$, then
$$x^2+7y^2\eq 14\da+r^2+7-14z^2\eq14\da-14(1-\da)\eq2\pmod 4$$
which is impossible. Thus $z\eq\da\pmod2$ and $x^2+7y^2\eq r^2+7\eq0\pmod8$.
Note that $x^2+7y^2\not=0$ since $7\nmid r$. Applying Lemma \ref{Lem3.2}(ii)
we can write $x^2+7y^2$ as $x_0^2+7y_0^2$ with $x_0$ and $y_0$ both odd.
Since $x_0^2\eq r^2\pmod7$, either $x_0$ or $-x_0$ has the form $14w+r$ with $w\in\Z$.
Write $y_0=2u+1$ and $z=2v+\da$ with $u,v\in\Z$. Then
$$56n+14\da+r^2+7\eq (14w+r)^2+7(2u+1)^2+14(2v+\da)^2$$
and hence
$n=p_3(u)+v(v+\da)+w(7w+r)/2$.
This proves the universality of $(7,r,2,2\da,1,1)$ over $\Z$.

The proof of Theorem 1.2 is now complete.
\qed

\section{Proof of Theorem \ref{Th1.3}}

For a positive definite integral ternary quadratic form $f(x,y,z)$ and an integer $n$, as usual we define
$$r(n,f):=|\{(x,y,z)\in\Z^3:\ f(x,y,z)=n\}|$$
and
\begin{equation*}r(n,\gen(f)):=\(\sum_{f^*\in \gen(f)}\frac{1}{|\Aut(f^*)|}\)^{-1}\sum_{f^*\in
  \gen(f)}\frac{r(n,f^*)}{|\Aut(f^*)|},
  \end{equation*}
  where the summation is over a set of representatives of the classes in $\gen(f)$,
  and $\Aut(f^*)$ is the group of integral isometries of $f^*$.

\begin{lemma}\label{Lem4.1} Let $f$ be a positive definite ternary quadratic form with determinant $d(f)$.
 Let $m\in\{1,2\}$ and suppose that $m$ is represented by the genus of $f$.
Then, for each prime $p\nmid 2md(f)$, we have
\begin{equation}\label{4.1} \frac{r(mp^2,\gen(f))}{r(m,\gen(f))}=p+1-\left(\frac{-md(f)}{p}\right),
 \end{equation}
 where $(\f{\cdot}p)$ denotes the Legendre symbol.
\end{lemma}
\begin{proof} By the Minkowski-Siegel formula \cite [pp.\,173-174]{Ki}, for any $n\in\Z^+$ we have
\begin{equation*} r(n,\gen(f))=2\pi\sqrt{\frac{n}{d(f)}}\ \prod_q\al_q(n,f),
\end{equation*}
where $q$ runs over all primes and $\al_q$ is the local density. As $p\nmid 2md(f)$, by \cite {Y} we have
\begin{align*}
\al_p(mp^2,f)& =1+\f{1}{p}-\f{1}{p^2}+\left(\frac{-md(f)}{p}\right)\f{1}{p^2},\\
\al_p(m,f)& =1+\left(\frac{-md(f)}{p}\right)\f{1}{p}.
\end{align*}
Thus
\begin{equation*}\frac{r(mp^2,\gen(f))}{r(m,\gen(f))}=p\frac{\al_p(mp^2,f)}{\al_p(m,f)}=p+1-\left(\frac{-md(f)}{p}\right).
\end{equation*}
This concludes the proof.
\end{proof}

\begin{lemma}\label{Lem4.2} Let $w=u^2+15v^2>0$ with $u,v\in\Z$ and $8\mid w$. Then $w=x^2+15y^2$ for some odd integers $x$ and $y$.
\end{lemma}
\begin{proof} Let $k$ be the $2$-adic order of $\gcd(u,v)$, and write $u=2^ku_0$ and $v=2^kv_0$ with $u_0,v_0\in\Z$ not all even.
If $k=0$, then both $u_0$ and $v_0$ are odd since $w$ is even. Below we assume $k>0$.

We observe the identity
\begin{equation*}4^2(x^2+15y^2)=(x-15y)^2+15(x+y)^2.
\end{equation*}
If $u_0\not\eq v_0\pmod2$, then $k\gs2$ (since $8\mid w$) and $4^2(u_0^2+15v_0^2)=s^2+15t^2$ with
$s=u_0-15v_0$ and $t=u_0+v_0$ both odd. For $j\in\N$, if $4^j(u_0^2+15v_0^2)=u_j^2+15v_j^2$
for some odd integers $u_j$ and $v_j$, then we may assume $u_j\eq v_j\pmod4$ without loss of generality
(otherwise we may replace $v_j$ by $-v_j$), and hence
$$4^{j+1}(u_0^2+15v_0^2)=4(u_j^2+15v_j^2)=u_{j+1}^2+15v_{j+1}^2$$
with $u_{j+1}=(u_j-15v_j)/2$ and $v_{j+1}=(u_j+v_j)/2$ both odd.
Thus, for some odd integers $u_k$ and $v_k$, we have
$$w=4^k(u_0^2+15v_0^2)=u_k^2+15v_k^2.$$
This concludes the proof.
\end{proof}

\medskip
\noindent{\it Proof of Theorem \ref{Th1.3}(i)}. (a) We first prove that $(7,7,3,1,1,1)$ is universal over $\Z$.
Let $n\in\N$. Clearly,
\begin{align*}&n=p_3(x)+7p_3(y)+\f{z(3z+1)}2
\\\iff&24n+25=3(2x+1)^2+21(2y+1)^2+(6z+1)^2.
\end{align*}
There are two classes in the genus of $x^2+3y^2+21z^2$ and the one not containing $x^2+3y^2+21z^2$ has the representative
\begin{equation}\label{4.2}\begin{aligned}x^2+6y^2+12z^2-6yz=&x^2+3\l(\frac{y}2-2z\r)^2+21\l(\frac{y}{2}\r)^2
\\=&x^2+3\l(\frac{y+3z}{2}\r)^2+21\l(\frac{y-z}{2}\r)^2.
\end{aligned}\end{equation}
If $24n+25=x^2+6y^2+12z^2-6yz$ with $x,y,z\in\Z$, then the equality modulo 4 yields
 $y(y-z)\equiv 0\pmod 2$. Thus, by (\ref{4.2}) and Lemma \ref{Lem3.1}, we have
 \begin{equation}\label{4.3}24n+25\in\{x^2+3y^2+21z^2:\ x,y,z\in\Z\}.
 \end{equation}

Now we claim that $24n+25=x^2+3y^2+21z^2$ for some $x,y,z\in\Z$ with $y^2+7z^2>0$.
This holds by (\ref{4.3}) if $24n+25$ is not a square.
Suppose that $24n+25=m^2$ with $m\in\Z^+$. Let $p$ be any  prime divisor of $m$. Clearly, $p\gs5$.
Note that $r(7^2,x^2+3y^2+21z^2)>2$ since $7^2=(\pm5)^2+3\times (\pm1)^2+21\times (\pm1)^2$.
If $p\not=7$ and $r(p^2,x^2+6y^2+12z^2-6yz)>2$, then $p^2=x^2+6y^2+12z^2-6yz$ for some $x,y,z\in\Z$ with $2\mid y(y-z)$ and $y^2+z^2>0$,
hence by (\ref{4.2}) we have $p^2=x^2+3u^2+21v^2$ for some $x,u,v\in\Z$ with $u^2+7v^2>0$, and thus $r(p^2,x^2+3y^2+21z^2)>2$.
By Lemma \ref{Lem4.1}, if $p\not=7$ then
\begin{align*}
\frac{r(p^2,\gen(x^2+3y^2+21z^2))}{r(1,\gen( x^2+3y^2+21z^2))}=p+1-\left(\frac{-7}{p}\right)
\end{align*}
and hence
\begin{equation*}r(p^2,x^2+3y^2+21z^2)+r(p^2,x^2+6y^2+12z^2-6yz)=4\l(p+1-\l(\frac{-7}{p}\r)\r)>4.
\end{equation*}
So we still have $r(p^2,x^2+3y^2+21z^2)>2$ if $r(p^2,x^2+6y^2+12z^2-6yz)\ls2$.
As $r(m^2,x^2+3y^2+21z^2)\gs r(p^2,x^2+3y^2+21z^2)>2$, we can write $24n+25=m^2$ as $x^2+3y^2+21z^2$ with $x,y,z\in\Z$ and $y^2+7z^2>0$.
This proves the claim.

By the claim, there are $x,y,z\in\Z$ such that $24n+25=x^2+3y^2+21z^2$ and $y^2+7z^2>0$. As $3y^2\not\eq 25\eq1\pmod4$, either $x$ or $z$ is odd.
If $2\nmid x$, then $3(y^2+7z^2)\eq 25-x^2\eq0\pmod8$ and hence by Lemma \ref{Lem3.2}(ii) we can write $y^2+7z^2$ as $y_0^2+7z_0^2$ with $y_0$ and $z_0$ both odd.
If $2\nmid z$, then $x^2+3y^2\eq 25-21z^2\eq4\pmod8$ and hence by Lemma \ref{Lem3.2}(i) we can write $x^2+3y^2$ as $x_1^2+3y_1^2$ with $x_1$ and $y_1$ both odd.
Thus $24n+25=a^2+3b^2+21c^2$ for some odd integers $a,b,c$. As $3\nmid a$, either $a$ or $-a$ has the form $6w+1$ with $w\in\Z$.
Write $b=2u+1$ and $c=2v+1$ with $u,v\in\Z$. Then
$$24n+25=(6w+1)^2+3(2u+1)^2+21(2v+1)^2$$
and hence $n=p_3(u)+7p_3(v)+w(3w+1)/2$. This proves the universality of $(7,7,3,1,1,1)$ over $\Z$.

(b) Let $n\in\N$ and $r\in\{1,3\}$. Clearly,
\begin{align*}&n=5p_3(x)+\f{y(3y+1)}2+\f{z(3z+r)}2
\\\iff&24n+r^2+16=15(2x+1)^2+(6y+1)^2+(6z+r)^2.
\end{align*}

There are two classes in the genus of $x^2+y^2+15z^2$, and the one not containing $x^2+y^2+15z^2$ has the representative
\begin{equation}\label{4.4}\begin{aligned}x^2+4y^2+4z^2-2yz=&x^2+\l(\f y2-2z\r)^2+15\l(\f y2\r)^2
 \\=&x^2+\l(2y-\frac{z}{2}\r)^2+15\l(\frac{z}{2}\r)^2.
\end{aligned}\end{equation}
If $24n+r^2+16=x^2+4y^2+4z^2-2yz$ with $x,y,z\in\Z$, then
$2\nmid x$ and $2\mid yz$. Thus, in view of (\ref{4.4}) and Lemma \ref{Lem3.1}, we have
\begin{equation}\label{4.5}24n+r^2+16\in\{x^2+y^2+15z^2:\ x,y,z\in\Z\}.
\end{equation}

We claim that $24n+r^2+16=x^2+y^2+15z^2$ for some $x,y,z\in\Z$ with $(x^2+15z^2)(y^2+15z^2)>0$.
This holds by (\ref{4.5}) if $24n+r^2+16$ is not a square.
Now suppose that $24n+r^2+16=m^2$ with $m\in\Z^+$. Let $p$ be any prime divisor of $m$.
Clearly, $p\gs5$. Note that $r(5^2,x^2+y^2+15z^2)>4$ since
$$5^2=(\pm5)^2+0^2+15\times0^2=0^2+(\pm 5)^2+15\times0^2=(\pm3)^2+(\pm4)^2+15\times 0^2.$$
If $r(p^2,x^2+4y^2+4z^2-2yz)>2$, then $p^2=x^2+4y^2+4z^2-2yz$ for some $x,y,z\in\Z$ with $2\mid yz$ and $y^2+z^2>0$,
hence by (\ref{4.4}) $p^2=x^2+u^2+15v^2$ for some $x,u,v\in\Z$ with $(x^2+15v^2)(u^2+15v^2)>0$, and thus $r(p^2,x^2+y^2+15z^2)>4$.
When $p>5$, by Lemma \ref{Lem4.1} we have
\begin{align*}
\frac{r(p^2, \gen(x^2+y^2+15z^2))}{r(1, \gen(x^2+y^2+15z^2))}=p+1-\l(\frac{-15}{p}\r)
\end{align*}
and hence
\begin{equation*}
r(p^2, x^2+y^2+15z^2)+2r(p^2,x^2+4y^2+4z^2-2yz)=8\l(p+1-\l(\frac{-15}{p}\r)\r)>50.
\end{equation*}
Thus we still have $r(p^2,x^2+y^2+15z^2)>4$ if $r(p^2,x^2+4y^2+4z^2-2yz)\ls2$.
As $r(m^2,x^2+y^2+15z^2)\gs r(p^2,x^2+y^2+15z^2)>4$, we can write $24n+r^2+16$ as $x^2+y^2+15z^2$ with $(x^2+15z^2)(y^2+15z^2)>0$.
This proves the claim.

By the claim, there are $x,y,z\in\Z$ such that $24n+r^2+16=x^2+y^2+15z^2$ and $(x^2+15z^2)(y^2+15z^2)>0$.
Since $15z^2\not\eq r^2\eq1\pmod 4$, either $x$ or $y$ is odd. Without any loss of generality, we assume that $2\nmid x$.
Since $y^2+15z^2>0$ and $y^2+15z^2\eq r^2-x^2\eq0\pmod8$, by Lemma \ref{Lem4.2} we can write $y^2+15z^2=y_0^2+15z_0^2$
with $y_0$ and $z_0$ both odd. Now, $24n+r^2+16=x^2+y_0^2+15z_0^2$. Since $x^2+y_0^2\eq r^2+1\pmod 3$, one of $x^2$ and $y_0^2$
is congruent to $r^2$ modulo 3 and the other one is congruent to $1$ modulo 3. Thus $x^2+y_0^2=(6u+r)^2+(6v+1)^2$ for some $u,v\in\Z$.
Write $z_0=2w+1$ with $v\in\Z$. Then
$$24n+r^2+16=(6u+r)^2+(6v+1)^2+15(2w+1)^2$$
and hence $n=u(3u+r)/2+v(3v+1)/2+5p_3(w)$. This proves the universality of $(5,5,3,r,3,1)$ over $\Z$.

(c) Let $n\in\N$. Apparently,
\begin{align*}&n=p_3(x)+5p_3(y)+z(3z+2)
\\\iff&24n+26=3(2x+1)^2+15(2y+1)^2+2(6z+2)^2.
\end{align*}

There are two classes in the genus of $2x^2+3y^2+15z^2$, and the one not containing $2x^2+3y^2+15z^2$ has the representative
\begin{equation}\label{4.6}g(x,y,z)=2x^2+5y^2+11z^2+2yz+2x(y-z)=2\l(x+v\r)^2+3(u-2v)^2+15u^2
\end{equation}
with $u=(y+z)/2$ and $v=(y-z)/2$.
If $24n+26=g(x,y,z)$ with $x,y,z\in\Z$, then $y\eq z\pmod2$, and hence by (\ref{4.6})
we have $24n+26=2a^2+3b^2+15c^2$ for some $a,b,c\in\Z$. So, in view of Lemma \ref{Lem3.1}, we always have
\begin{equation}\label{4.7}24n+26\in\{2x^2+3y^2+15z^2:\ x,y,z\in\Z\}.
\end{equation}

We claim that $24n+26=2x^2+3y^2+15z^2$ for some $x,y,z\in\Z$ with $y^2+5z^2>0$.
This holds by (\ref{4.7}) if $12n+13$ is not a square.
Now suppose that $12n+13=m^2$ with $m\in\Z^+$. Let $p$ be any prime divisor of $m$.
Clearly, $p\gs5$. Note that $r(2\times5^2,2x^2+3y^2+15z^2)>2$ since
$$2\times5^2=2\times (\pm5)^2+3\times0^2+15\times0^2=2(\pm1)^2+3(\pm4)^2+30\times0^2.$$
If $r(2p^2,g(x,y,z))>2$, then $2p^2=g(x,y,z)$ for some $x,y,z\in\Z$ with $y^2+z^2>0$,
hence by (\ref{4.6}) $2p^2=2x^2+3b^2+15c^2$ for some $x,b,c\in\Z$ with $b^2+c^2>0$, and thus $r(2p^2,2x^2+3y^2+15z^2)>2$.
When $p>5$, by Lemma \ref{Lem4.1} we have
\begin{align*}
\frac{r(2p^2, \gen(2x^2+3y^2+15z^2))}{r(2, \gen(2x^2+3y^2+15z^2))}=p+1-\l(\frac{-5}{p}\r)
\end{align*}
and hence
\begin{align*}
r(2p^2, 2x^2+3y^2+15z^2)+2r(2p^2,g(x,y,z)) =6\l(p+1-\l(\frac{-5}{p}\r)\r)>40.
\end{align*}
Thus we still have $r(2p^2,2x^2+3y^2+15z^2)>2$ if $r(2p^2,g(x,y,z))\ls2$.
As $r(2m^2,2x^2+3y^2+15z^2)\gs r(2p^2,2x^2+3y^2+15z^2)>2$, we can write $24n+26$ as $2x^2+3y^2+15z^2$ with $y^2+5z^2>0$.
This proves the claim.

By the claim, there are $x,y,z\in\Z$ such that $24n+26=2x^2+3(y^2+5z^2)$ and $y^2+5z^2>0$. By \cite[Lemma 2.1]{S15}, $y^2+5z^2=y_0^2+5z_0^2$
for some integers $y_0$ and $z_0$ not all divisible by $3$. Without any loss of generality, we simply assume that $3\nmid y$ or $3\nmid z$.
Note that $3\nmid x$ and $y\eq z\pmod2$. If $3\nmid yz$, then $\ve_1 y\eq \ve_2 z\eq x\pmod3$ for some $\ve_1,\ve_2\in\{\pm1\}$.
If $3\mid y$ and $3\nmid z$ then $x+y+\ve z\eq0\pmod3$ for some $\ve\in\{\pm1\}$; similarly, if $3\nmid y$ and $3\mid z$ then
$x+\ve y+z\eq0\pmod3$. So, without loss of generality we may suppose that $x+y+z\eq0\pmod3$ (otherwise we adjust signs of $x,y,z$ suitably to meet our purpose).
If $y\eq z\eq0\pmod2$, then $2x^2\eq 26\pmod4$, hence $2\nmid x$ and $y\eq z\pmod 4$ since $y^2+5z^2\eq0\pmod8$, therefore
\begin{equation}\label{4.8}2x^2+3y^2+15z^2=2\l(\f{y-5z}2\r)^2+3\l(\f{2x+5y+5z}6\r)^2+15\l(\f{2x-y-z}6\r)^2
\end{equation}
with $(2x+5y+5z)/6$ and $(2x-y-z)/6$ both odd.

By the above, $24n+26=2a^2+3b^2+15c^2$ for some $a,b,c\in\Z$ with $2\nmid bc$. As $3\nmid a$ and $2a^2\eq 26-3-15\eq0\pmod8$,
$a$ or $-a$ has the form $2(3w+1)$ with $w\in\Z$. Write $b=2u+1$ and $c=2v+1$ with $u,v\in\Z$. Then
$$24n+26=2(2(3w+1))^2+3(2u+1)^2+15(2v+1)^2$$
and so $n=p_3(u)+5p_3(v)+w(3w+2)$. This proves the universality of $(6,4,5,5,1,1)$ over $\Z$.
\qed

\medskip
\noindent{\it Proof of Theorem \ref{Th1.3}(ii)}. (a) Let $n\in\N$ and $r\in\{1,2\}$. It is easy to see that
\begin{align*}&n=p_3(x)+5\f{y(3y+1)}2+z(3z+r)
\\\iff&24n+2r^2+8=3(2x+1)^2+5(6y+1)^2+2(6z+r)^2.
\end{align*}

As mentioned in part (b) of the proof of Theorem 1.3(i), there are two classes in the genus of $x^2+y^2+15z^2$, and the one not containing $x^2+y^2+15z^2$ has the representative
$x^2+4y^2+4z^2-2yz$. If $12n+r^2+4=x^2+4y^2+4z^2-2yz$ with $x,y,z\in\Z$, then $2\mid yz$ since $r^2\not\eq x^2-2\pmod4$.
Thus, in view of (\ref{4.4}) and Lemma \ref{Lem3.1}, $12n+r^2+4=x^2+y^2+15z^2$ for some $x,y,z\in\Z$.
If $x\eq y\pmod2$, then $z\eq r\pmod2$, $x^2+y^2\eq r^2-15z^2\eq 2r^2\pmod4$ and hence $x\eq y\eq r\eq z\pmod2$.
So, $x$ or $y$ has the same parity with $z$. Without loss of generality we may assume that $y\eq z\pmod 2$.
Since $y^2+15z^2\eq0\pmod4$, we have $x\eq r\pmod2$. If $r=2$ and $y^2+15z^2=0$,  then $12n+r^2+4=0^2+x^2+15\times0^2$
with $x\eq 0\eq r\pmod2$ and $x^2+15\times0^2>0$. If $r=1$, then $12n^2+r^2+4=12n+5$ is congruent to 2 modulo 3 and hence not a square.
Thus, without loss of generality we may assume that $y^2+15z^2>0$.

Observe that
$$24n+2r^2+8=2(x^2+y^2+15z^2)=2x^2+3u^2+5v^2$$
with $u=(y+5z)/2$ and $v=(y-3z)/2$ both odd. Since $3u^2+5v^2\eq 2r^2-2x^2\eq0\pmod8$ and $2(3u^2+5v^2)=y^2+15z^2>0$,
by Lemma \ref{Lem3.2}(iii) we can write $3u^2+5v^2$ as $3y_0^2+5z_0^2$ with $y_0$ and $z_0$ both odd.
As $2(x^2+z_0^2)\eq 2x^2+5z_0^2\eq 2r^2+8\pmod3$, we have $x^2+z_0^2\eq r^2+1\eq2\pmod3$ and hence
we may write $x$ or $-x$ as $6u+r$, $z_0$ or $-z_0$ as $6v+1$, and $y_0=2w+1$, where $u,v,w$ are integers.
Therefore
$$24n+2r^2+8=2x^2+3y_0^2+5z_0^2=2(6u+r)^2+3(2w+1)^2+5(6v+1)^2$$
and hence $n=u(3u+r)/2+5v(3v+1)/2+p_3(w)$. This proves the universality of $(15,5,6,2r,1,1)$ over $\Z$.

(b) Let $n\in\N$, $s\in\{1,3,5\}$ and $t\in\{1,2\}$ with $(s,t)\not=(5,2)$. Obviously,
\begin{align*}&n=p_3(x)+\f{y(5y+s)}2+z(3z+t)
\\\iff&120n+3s^2+10t^2+15=15(2x+1)^2+3(10y+s)^2+10(6z+t)^2.
\end{align*}

There are two classes in the genus of $3x^2+10y^2+15z^2$, and the one not containing $3x^2+10y^2+15z^2$ has the representative
\begin{equation}\label{4.9}\begin{aligned}g(x,y,z)=&7x^2+7y^2+12z^2+6(x+y)z+4xy
\\=&3\l(\f{x+y}2+2z\r)^2+10\l(\f{x-y}2\r)^2+15\l(\f{x+y}2\r)^2.
\end{aligned}\end{equation}
If $120n+3s^2+10t^2+15=g(x,y,z)$ with $x,y,z\in\Z$, then we obviously have $x\eq y\pmod2$.
Thus, in view of (\ref{4.9}) and Lemma \ref{Lem3.1}, $120n+3s^2+10t^2+15=3x^2+10y^2+15z^2$ for some $x,y,z\in\Z$.
If $x=z=0$, then $120n+3s^2+10t^2+15=10y^2$, hence $(s,t)=(5,1)$ and $y^2=12n+10\eq2\pmod 4$ which is impossible.
So $x^2+5z^2>0$, and hence by \cite[Lemma 2.1]{S15} we can rewrite $x^2+5z^2$ as $x_0^2+5z_0^2$ with $x_0,z_0\in\Z$ not all divisible by 3.
Without loss of generality, we simply assume that $3\nmid x$ or $3\nmid z$. Note that $3\nmid y$ since $3\nmid t$.
If $3\nmid xz$, then $\ve_1 x\eq y\eq\ve_2 z$ for some $\ve_1,\ve_2\in\{\pm1\}$.
If $3\mid x$ and $3\nmid z$, then $x+y+\ve z\eq0\pmod3$ for some $\ve\in\{\pm1\}$.
If $3\nmid x$ and $3\mid z$, then $\ve x+y+z\eq0\pmod3$ for some $\ve\in\{\pm1\}$.
Without loss of generality, we just assume that $x+y+z\eq0\pmod3$ (otherwise we may adjust signs of $x,y,z$ suitably).
Note that $x\eq z\pmod2$ and we have the identity
\begin{equation}\label{4.10}3x_1^2+10y_1^2+15z_1^2=3x^2+10y^2+15z^2,
\end{equation}
where
$$x_1=\f{x+10y-5z}6,\ \ y_1=\f{x+z}2\ \ \t{and}\ \ z_1=\f{x-2y-5z}6$$ are all integral.

If $x\eq z\eq1\pmod2$, then $10y^2=120n+3s^2+10t^2+15-3x^2-15z^2\eq 10t^2\pmod4$ and hence $y\eq t\pmod2$.

Now suppose that $x\eq z\eq0\pmod2$. Then $2y^2\eq 10y^2\eq 3s^2+10t^2+15\eq 2(t^2+1)\pmod4$ and hence $y\not\eq t\pmod2$.
Observe that
$$2t^2+2\eq 120n+3s^2+10t^2+15=3x^2+10y^2+15z^2\eq x^2+z^2+2(t+1)^2\pmod8$$
and hence
$$y_1=\f{x+z}2\eq\l(\f x2\r)^2+\l(\f z2\r)^2=\f{x^2+z^2}4\eq t\pmod2.$$
Thus
$$z_1=x_1-2y\eq x_1\eq\f{x+z}2-3z+5y\eq t+y\eq1\pmod2.$$

In view of the above, there are integers $x,y,z\in\Z$ with $x\eq z\eq1\pmod2$ and $y\eq t\pmod2$ such that $120n+3s^2+10t^2+15=3x^2+10y^2+15z^2$.
 Clearly, $y$ or $-y$ has the form $6v+t$ with $v\in\Z$. Write $z=2w+1$ with $w\in\Z$. Since $x^2\eq s^2\pmod5$, we can write $x$ or $-x$ as $10u+s$ with $w\in\Z$.
 Therefore
 $$120n+3s^2+10t^2+15=3(10u+s)^2+10(6v+t)^2+15(2w+1)^2$$
 and hence $n=p_3(w)+u(5u+s)/2+v(3v+t)$. This proves the universality of $(6,2t,5,s,1,1)$ over $\Z$.
 \qed

\section{Proof of Theorem \ref{Th1.4}}

B. W. Jones and G. Pall \cite{Jones} proved the following celebrated result.
\begin{lemma}\label{Lem5.1} Let $n\in \N$ with $8n+1$ not a square. Then
\begin{align*}&|\{(x,y,z)\in\Z^3:\ x^2+y^2+z^2=8n+1\ \&\ 4\mid x\}|
\\=&|\{(x,y,z)\in\Z^3:\ x^2+y^2+z^2=8n+1\ \&\ x\eq2\pmod4\}| >0.
\end{align*}
\end{lemma}

A. G. Earnest \cite {E80,E84} showed the following useful result.
\begin{lemma}\label{Lem5.2} Let $c$ be a primitive spinor exceptional integer for the genus of a positive ternary quadratic form $f(x,y,z)$, and let $S$ be a spinor genus containing $f$.
Let $s$ be a fixed positive integer relatively prime to $2d(f)$ for which $cs^2$ can be primitively represented by $S$.
If $t\in\Z^+$ is relatively prime to $2d(f)$, then $ct^2$ can be primitively represented by $S$ if and only if
$$\left(\f{-c d(f)}{s}\right)= \left(\f{-cd(f)}{t}\right).$$
\end{lemma}

\medskip
\noindent{\it Proof of Theorem \ref{Th1.4}}. Fix $n\in\N$ . Clearly,
\begin{align*}n=p_3(x)+y^2+2z(4z+1)\iff 8n+2=(2x+1)^2+8y^2+(8z+1)^2.
\end{align*}
So, it suffices to show that $8n+2=x^2+y^2+8z^2$ for some $x,y,z\in\Z$ with $x\eq\pm1\pmod8$.

\vspace*{4pt}\noindent{\bf Case 1.} $n$ is not twice a triangular number.

In this case, $4n+1$ is not a square. If $2\mid n$, then by  Lemma \ref{Lem5.1} we can write $4n+1$ as $x^2+y^2+z^2$
with $2\nmid x$, $2\mid y$ and $z\eq 2\pmod4$. If $2\nmid n$, then there are $x,y,z\in\Z$ with $2\nmid x$ and $y\eq z\eq0\pmod2$
such that $4n+1=x^2+y^2+z^2$ and hence $y\not\eq z\pmod4$ since $y^2+z^2\eq 5-x^2\eq4\pmod8$. So we can always write $4n+1=x^2+y^2+z^2$
with $2\nmid x$, $2\mid y$ and $z\eq 2n-2\pmod4$, hence
$$8n+2=2(x^2+y^2+z^2)=(x+y)^2+(x-y)^2+8\l(\f z2\r)^2$$
with $z/2\eq n-1\pmod2$. Thus
$$(x+y)^2+(x-y)^2\eq 8n+2-8(n-1)=10\not\eq 3^2+3^2\pmod{16}$$
and hence $x+\ve y\eq\pm1\pmod8$ for some $\ve\in\{\pm1\}$.

\vspace*{4pt}\noindent{\bf Case 2.} $n=2p_3(m)$ with $m\in\N$, and $2m+1$ has no prime factor of the form $4k+3$.

In this case, $2m+1$ can be expressed as the sum of two squares. If $4\mid m$, then
$$8n+2=2(2m+1)^2=(2m+1)^2+(2m+1)^2+8\times0^2$$
with $2m+1\eq1\pmod8$. If $4\nmid m$, then $2m+1=u^2+(2v)^2$ for some odd integers $u$ and $v$, and hence
\begin{align*}8n+2=&2(u^2+4v^2)^2=2((u^2-4v^2)^2+(4uv)^2)
\\=&(u^2-4v^2+4uv)^2+(u^2-4v^2-4uv)^2+8\times0^2
\end{align*}
with $u^2-4v^2\pm4uv\eq 1\pmod8$.

\vspace*{4pt}\noindent{\bf Case 3}. $n=2p_3(m)$ with $m\in\N$, and $2m+1$ has a prime factor $p\eq3\pmod4$.

By Lagrange's four-square theorem, we can write $p=a^2+b^2+c^2+d^2$, where $a$ is an even number and $b,c,d$ are odd numbers.
Thus
\begin{align*}p^2=&(a^2+b^2-c^2-d^2)^2+4(a^2+b^2)(c^2+d^2)
\\=&(a^2+b^2-c^2-d^2)^2+(2ac+2bd)^2+(2ad-2bc)^2
\end{align*}
and hence $(2m+1)^2=x^2+(2y)^2+(2z)^2$ for some odd integers $x,y,z$. Observe that
$$8n+2=2(2m+1)^2=(x+2y)^2+(x-2y)^2+8z^2$$
and $(x+2y)^2+(x-2y)^2\eq 2-8z^2\eq10\not\eq3^2+3^2\pmod{16}$.
So one of $x+2y$ and $x-2y$ is congruent to $1$ or $-1$ modulo $8$.

Now we give an alternative approach to Case 3. There are three classes in the genus of $x^2+y^2+32z^2$ with the three representatives
\begin{align*}f_1(x,y,z)=&x^2+y^2+32z^2,
\\ f_2(x,y,z)=&2x^2+2y^2+9z^2+2yz-2zx,
\\ f_3(x,y,z)=&x^2+4y^2+9z^2-4yz.
\end{align*}
The class of $f_1$ and the class of $f_2$ constitute a spinor genus while another spinor genus in the genus only contains the class of $f_3$.
Since $2$ is a primitive spinor exceptional integer for this genus, by Lemma \ref{Lem5.2} we can write $2p^2$ as
$$f_3(u,v,w)=u^2+4v^2+9w^2-4vw=u^2+(2v-w)^2+8w^2$$ with $u,v,w\in \Z$. Since $2\nmid uw$,  we see that $8n+2=2(2m+1)^2=a^2+b^2+8c^2$ for some odd integers $a,b,c$.
As $a^2+b^2\eq2-8c^2\eq10\not\eq3^2+3^2\pmod{16}$, $a$ or $b$ is congruent to $1$ or $-1$ modulo $8$.
This concludes our discussion of Case 3.

In view of the above, we have completed the proof of Theorem \ref{Th1.4}.
\qed



\begin{thebibliography}{99}
\bibitem{B}  [10.1090/stml/034]
\newblock B. C. Berndt,
\newblock \emph{Number Theory in the Spirit of Ramanujan},
\newblock Amer. Math. Soc., Providence, RI, 2006. 

\bibitem{C}
\newblock J. W. S. Cassels,
\newblock \emph{Rational Quadratic Forms},
\newblock Academic Press, London, 1978. 

\bibitem{D27}  [10.2307/2370770]
\newblock L. E. Dickson,
\newblock {Quaternary quadratic forms representing all integers},
\newblock \emph{Amer. J. Math}., \textbf{49} (1927), 39--56. 

\bibitem{D39}
\newblock L. E. Dickson,
\newblock \emph{Modern Elementary Theory of Numbers},
\newblock Univ. of Chicago Press, Chicago, 1939. 

\bibitem{E80}  [10.2140/pjm.1980.90.325]
\newblock A. G. Earnest,
\newblock {Congruence conditions on integers represented by ternary quadratic forms},
\newblock \emph{Pacific J. Math}., \textbf{90} (1980), 325--333. 

\bibitem{E84}  [10.1017/S0027763000020717]
\newblock A. G. Earnest,
\newblock {Representation of spinor exceptional integers by ternary quadratic forms},
\newblock \emph{Nagoya Math. J}., \textbf{93} (1984), 27--38. 

\bibitem{GS}
\newblock F. Ge and Z.-W. Sun,
\newblock On some universal sums of generalized polygonals,
\newblock \emph{Colloq. Math}., \textbf{145} (2016), 149--155. 

\bibitem{GPS}
\newblock S. Guo, H. Pan and Z.-W. Sun,
\newblock Mixed sums of squares and triangular numbers (II),
\newblock \emph{Integers}, \textbf{7} (2007), A56, 5pp (electronic). 

\bibitem{Jagy96}  [10.4064/aa-77-4-361-367]
\newblock W. C. Jagy,
\newblock {Five regular or nearly-regular ternary quadratic forms},
\newblock \emph{Acta Arith}., \textbf{77} (1996), 361--367. 

\bibitem{JKS}  [10.1112/S002557930001264X]
\newblock W. C. Jagy, I. Kaplansky and A. Schiemann,
\newblock {There are 913 regular ternary forms},
\newblock \emph{Mathematika}, \textbf{44} (1997), 332--341. 

\bibitem{Jagy14}
\newblock W. C. Jagy,
\newblock \emph{Integral Positive Ternary Quadratic Forms},
\newblock Lecture Notes, 2014. Available from: \url{http://zakuski.math.utsa.edu/~kap/Jagy_Encyclopedia.pdf}.

\bibitem{Jones}  [10.1007/BF02547347]
\newblock B. W. Jones and G. Pall,
\newblock {Regular and semi-regular positive ternary quadratic forms},
\newblock \emph{Acta Math}., \textbf{70} (1939), 165--191. 

\bibitem{JOS}  [10.1142/S1793042119500350]
\newblock J. Ju, B.-K. Oh and B. Seo,
\newblock {Ternary universal sums of generalized polygonal numbers},
\newblock \emph{Int. J. Number Theory}, \textbf{15} (2019), 655--675. 

\bibitem{Ki}  [10.1017/CBO9780511666155]
\newblock Y. Kitaoka,
\newblock \emph{Arithmetic of Quadratic Forms},
\newblock Cambridge Tracts in Math., Vol. 106, Cambridge, 1993. 

\bibitem{Oh}  [10.4134/JKMS.2011.48.4.837]
\newblock B.-K. Oh,
\newblock {Ternary universal sums of generalized pentagonal numbers},
\newblock \emph{J. Korean Math. Soc}., \textbf{48} (2011), 837--847. 

\bibitem{OS}  [10.1016/j.jnt.2008.10.002]
\newblock B.-K. Oh and Z.-W. Sun,
\newblock {Mixed sums of squares and triangular numbers (III)},
\newblock \emph{J. Number Theory}, \textbf{129} (2009), 964--969. 

\bibitem{Oto}
\newblock O. T. O'Meara,
\newblock \emph{Introduction to Quadratic Forms},
\newblock Springer, New York, 1963. 

\bibitem{R}
\newblock S. Ramanujan,
\newblock On the expression of a number in the form $ax^2+by^2+cz^2+dw^2$,
\newblock \emph{Proc. Cambridge Philos. Soc}., \textbf{19} (1917), 11--21.

\bibitem{S07}  [10.4064/aa127-2-1]
\newblock Z.-W. Sun,
\newblock {Mixed sums of squares and triangular numbers},
\newblock \emph{Acta Arith}., \textbf{127} (2007), 103--113. 

\bibitem{S15}  [10.1007/s11425-015-4994-4]
\newblock Z.-W. Sun,
\newblock {On universal sums of polygonal numbers},
\newblock \emph{Sci. China Math}., \textbf{58} (2015), 1367--1396. 

\bibitem{S16}  [10.1016/j.jnt.2015.10.014]
\newblock Z.-W. Sun,
\newblock {A result similar to Lagrange's theorem},
\newblock \emph{J. Number Theory}, \textbf{162} (2016), 190--211. 

\bibitem{S17}  [10.1016/j.jnt.2016.07.024]
\newblock Z.-W. Sun,
\newblock {On $x(ax+1)+y(by+1)+z(cz+1)$ and $x(ax+b)+y(ay+c)+z(az+d)$},
\newblock \emph{J. Number Theory}, \textbf{171} (2017), 275--283. 


\bibitem{S17b}
\newblock Z.-W. Sun,
\newblock Sequence A286944 in OEIS, 2017.
\newblock Available from: \url{http://oeis.org/A286944}.

\bibitem{S17a} [10.1007/s11425-017-9354-4]
\newblock Z.-W. Sun,
\newblock {Universal sums of three quadratic polynomials},
\newblock \emph{Sci. China Math.}, 2018.
 Available from: \url{https://doi.org/10.1007/s11425-017-9354-4}.
 See also {\tt arXiv:1502.03056}.

\bibitem{S18}
\newblock Z.-W. Sun,
\newblock On universal sums $x(ax+b)/2+y(cy+d)/2+z(ez+f)/2$,
\newblock \emph{Nanjing Univ. J. Math. Biquarterly}, \textbf{35} (2018), 85--199.

\bibitem{Y}  [10.1006/jnth.1998.2258]
\newblock T. Yang,
\newblock {An explicit formula for local densities of quadratic forms},
\newblock \emph{J. Number Theory}, \textbf{72} (1998), 309--356. 

\end{thebibliography}
\end{document}